\documentclass[11pt, final]{article}
\usepackage{a4}
\usepackage{amsmath}%
\usepackage{amstext}%
\usepackage{amssymb}%
\usepackage{showkeys}%
\usepackage{epsfig}%
\usepackage{cite}
\usepackage{tikz}
\usepackage{subfig}
\usepackage{caption}
\usepackage{multirow}
\usepackage{booktabs}
\usepackage{floatrow}

\usepackage{listings}
\usepackage{geometry}
\newtheorem{theorem}{Theorem}

\newtheorem{axiom}{Axiom}

\newtheorem{definition}[axiom]{Definition}

\newtheorem{lemma}[theorem]{Lemma}

\newtheorem{proposition}[theorem]{Proposition}
\newtheorem{rem}[theorem]{Remark}
\newenvironment{remark}{\rem\rm}{\endrem}

\newcounter{unnumber}

\newtheorem{prf}[unnumber]{Proof.}
\newenvironment{proof}{\prf\rm}{\hfill{$\blacksquare$}\endprf}
\newcommand{\R}{\mathbb{R}}%
\newcommand{\N}{\mathbb{N}}%
\newcommand{\ol}{\overline}%
\DeclareMathOperator*\inte{int}%
\DeclareMathOperator*\sqri{sqri}%
\DeclareMathOperator*\ri{ri}%
\DeclareMathOperator*\dom{dom}%
\DeclareMathOperator*\ran{ran}%
\DeclareMathOperator*\id{Id}%
\DeclareMathOperator*\argmin{argmin}

\title{Levenberg-Marquardt dynamics associated to variational inequalities}

\author{Radu Ioan Bo\c{t} \thanks{University of Vienna, Faculty of Mathematics, Oskar-Morgenstern-Platz 1, A-1090 Vienna, Austria,
email: radu.bot@univie.ac.at.} \and
Ern\"{o} Robert Csetnek \thanks {University of Vienna, Faculty of Mathematics, Oskar-Morgenstern-Platz 1, A-1090 Vienna, Austria,
email: ernoe.robert.csetnek@univie.ac.at. Research supported by FWF (Austrian Science Fund), Lise Meitner Programme, project M 1682-N25.}}

\begin{document}
\maketitle

\noindent \textbf{Abstract.} In connection with the optimization problem $$\inf_{x\in\argmin\Psi}\{\Phi(x)+\Theta(x)\},$$ 
where $\Phi$ is a proper, convex and lower semicontinuous function and $\Theta$ and $\Psi$ are convex and smooth functions defined on a real Hilbert space, 
we investigate the asymptotic behavior of the trajectories of the nonautonomous 
Levenberg-Marquardt dynamical system 
\begin{equation*}\left\{
\begin{array}{ll}
v(t)\in\partial\Phi(x(t))\\
\lambda(t)\dot x(t) + \dot v(t) + v(t) + \nabla \Theta(x(t))+\beta(t)\nabla \Psi(x(t))=0,
\end{array}\right.\end{equation*}
where $\lambda$ and $\beta$ are functions of time controlling the velocity and  the penalty term, respectively. We show weak convergence of the generated trajectory to an optimal solution  
as well as convergence of the objective function values along the trajectories, provided $\lambda$ is monotonically decreasing, $\beta$ satisfies a growth condition and a relation expressed 
via the Fenchel conjugate of $\Psi$ is fulfilled. When the objective function is assumed to be strongly convex, we can even show 
strong convergence of the trajectories.
\vspace{1ex}

\noindent \textbf{Key Words.} nonautonomous systems, Levenberg-Marquardt dynamics, regularized Newton-like dynamics, 
Lyapunov analysis, convex optimization, variational inequalities, penalization techniques\vspace{1ex}

\noindent \textbf{AMS subject classification.} 34G25, 47J25, 47H05, 90C25

\section{Introduction}\label{sec1}

Throughout this manuscript ${\cal H}$ is assumed to be a real Hilbert space endowed with inner product $\langle\cdot,\cdot\rangle$
and associated norm $\|\cdot\|\!=\!\sqrt{\langle \cdot,\cdot\rangle}$. When $T: {\cal H}\rightarrow {\cal H}$ is a $C^1$ operator 
with derivative $T'$, the solving of the equation 
$$\mbox{find }x\in {\cal H} \mbox{ such that } Tx=0$$ can be approached by the classical Newton method, 
which generates an approximating sequence $(x_n)_{n\geq 0}$ of a solution of the operator equation through
$$T(x_n)+T'(x_n)(x_{n+1}-x_n)=0 \ \forall n \geq 0.$$
In order to overcome the fact that the classical Newton method assumes the solving of an equation which
is in general not well-posed, one can use instead the Levenberg-Marquardt method 
 $$T(x_n)+\Big(\lambda_n\id+T'(x_n)\Big)\left(\frac{x_{n+1}-x_n}{\Delta t_n}\right)=0 \ \forall n \geq 0,$$
where $\id : {\cal H} \rightarrow {\cal H}$ denotes the identity operator on ${\cal H}$, $\lambda_n$ a regularizing parameter and $\Delta t_n > 0$ the step size.

When $T:{\cal H}\rightrightarrows {\cal H}$ is a (set-valued) maximally monotone operator, 
Attouch and Svaiter showed in \cite{att-sv2011} that the above Levenberg-Marquardt algorithm can be seen as a time discretization of the dynamical system 
\begin{equation}\label{att-sv-dyn}\left\{
\begin{array}{ll}
v(t)\in T(x(t))\\
\lambda(t)\dot x(t) + \dot v(t) + v(t)=0
\end{array}\right.\end{equation}
for approaching the inclusion problem 
\begin{equation}\label{o-tx}\mbox{find }x\in {\cal H} \mbox{ such that } 0\in Tx.\end{equation}
This includes as a special instance the problem of minimizing a proper, convex and lower semicontinuous function, when  $T$ is taken as its convex subdifferential. Later on, this investigation has been continued in \cite{abbas-att-sv} in the context of minimizing the  sum of a proper, convex and lower semicontinuous function with a convex and smooth one. 

In the spirit of \cite{att-sv2011}, we approach in this paper the optimization problem  
\begin{equation}\label{opt-intr}\inf_{x\in\argmin\Psi}\{\Phi(x)+\Theta(x)\},\end{equation}
where $\Phi:{\cal H}\rightarrow \R\cup\{+\infty\}$ is a proper, convex and lower semicontinuous function 
and $\Psi,\Theta:{\cal H}\rightarrow \R$ are convex and smooth functions, via the 
Levenberg-Marquardt dynamical system 
\begin{equation}\label{dyn-sys-intr}\left\{
\begin{array}{ll}
v(t)\in\partial\Phi(x(t))\\
\lambda(t)\dot x(t) + \dot v(t) + v(t) + \nabla \Theta(x(t))+\beta(t)\nabla \Psi(x(t))=0,
\end{array}\right.\end{equation} 
where $\lambda$  and $\beta$ are functions of time controlling the velocity and the penalty term, respectively. 

If $\partial \Phi + N_{\argmin\Psi}$ is maximally monotone, then determining an optimal solution $x \in {\cal H}$ of 
\eqref{opt-intr} means nothing else than solving the subdifferential inclusion problem
\begin{equation}\label{moninclusion}
\mbox{find} \ x \in {\cal H} \ \mbox{such that} \ 0 \in \partial \Phi(x) + \nabla\Theta(x)+ N_{\argmin\Psi}(x)
\end{equation}
or, equivalently, solving the variational inequality
\begin{equation}\label{vi}
\mbox{find} \ x \in {\argmin\Psi}  \ \mbox{and} \ v \in \partial \Phi(x) \ \mbox{such that} 
\ \langle v+\nabla\Theta(x), y-x \rangle \geq 0 \ \forall y \in \argmin\Psi.
\end{equation}

We show weak convergence of the trajectory $x(\cdot)$ generated by \eqref{dyn-sys-intr} to an optimal solution of \eqref{opt-intr} as well as 
convergence of the objective function values along the trajectory to the optimal objective value, 
provided the assumption \begin{equation}\label{ass-att-cz}\forall p\in\ran N_{\argmin\Psi} \ \int_0^{+\infty} 
\beta(t)\left[\Psi^*\left(\frac{p}{\beta(t)}\right)-\sigma_{\argmin\Psi}\left(\frac{p}{\beta(t)}\right)\right]dt<+\infty
\end{equation} is fulfilled and the functions $\lambda, \beta$ satisfy some mild conditions. If the objective 
function of \eqref{opt-intr} is strongly convex, the trajectory $x(\cdot)$ converges even strongly to the unique optimal solution 
of \eqref{opt-intr}. 

The condition \eqref{ass-att-cz} has its origins in the paper of Attouch and Czarnecki \cite{att-cza-10}, where the solving of 
\begin{equation}\inf_{x\in\argmin\Psi} \Phi(x),\end{equation}
for $\Phi, \Psi:{\cal H}\rightarrow \R\cup\{+\infty\}$ proper, convex and lower semicontinuous functions, 
is approached through the nonautonomous first order dynamical system 
\begin{equation}\label{1ord-att-cz}0\in\dot x(t)+\partial \Phi(x(t))+\beta(t)\partial \Psi(x(t)),\end{equation}
by assuming that the penalizing function $\beta:[0,+\infty)\rightarrow(0,+\infty)$  tend to $+\infty$ as $t\rightarrow +\infty$. 
Several ergodic and nonergodic convergence results have been reported in \cite{att-cza-10} under the key assumption \eqref{ass-att-cz}.

The paper of Attouch and Czarnecki \cite{att-cza-10} was the starting point of a remarkable number of research articles devoted 
to penalization techniques for solving optimization problems of type \eqref{opt-intr}, but also generalizations of the latter in form of variational inequalities expressed with maximal monotone operators (see \cite{att-cza-10, att-cza-peyp-c, att-cza-peyp-p, noun-peyp, peyp-12, b-c-penalty-svva, 
b-c-penalty-vjm, banert-bot-pen, b-c-dyn-pen, att-cab-cz, att-maing, b-c-dyn-sec-ord-pen}). In the literature enumerated above, 
the monotone inclusions problems have been approached either through continuous dynamical systems or through their 
discrete counterparts formulated as splitting algorithms. We speak in both cases about methods of penalty type, which means in this context that the operator describing the underlying set of the  variational inequality under investigation is evaluated as a penalty functional. In the above-listed references one can find more general formulations of the key assumption \eqref{ass-att-cz}, but also further examples for which these conditions are satisfied. In Remark \ref{phi,psi0} and Remark \ref{discrete} we provide more insights into the relations of the dynamical system \eqref{dyn-syst} to other continuous systems (and their discrete counterparts) from the literature.

\section{Preliminaries}\label{sec2}

In this section we present some preliminary definitions, results and tools that will be useful throughout the paper. 
We consider the following definition of an absolutely continuous function. 

\begin{definition}\label{abs-cont} \rm (see, for instance, \cite{att-sv2011, abbas-att-sv}) A function $x:[0,b]\rightarrow {\cal H}$ 
(where $b>0$) is said to be absolutely continuous if one of the 
following equivalent properties holds:
\begin{enumerate}

\item[(i)] there exists an integrable function $y:[0,b]\rightarrow {\cal H}$ such that $$x(t)=x(0)+\int_0^t y(s)ds \ \ \forall t\in[0,b];$$

\item[(ii)] $x$ is continuous and its distributional derivative is Lebesgue integrable on $[0,b]$; 

\item[(iii)] for every $\varepsilon > 0$, there exists $\eta >0$ such that for any finite family of intervals $I_k=(a_k,b_k) \subseteq [0,b]$ 
we have the implication
$$\left(I_k\cap I_j=\emptyset \mbox{ and }\sum_k|b_k-a_k| < \eta\right)\Longrightarrow \sum_k\|x(b_k)-x(a_k)\| < \varepsilon.$$
\end{enumerate}
A function $x:[0,+\infty) \rightarrow {\cal H}$ (where $b>0$) is said to be locally absolutely continuous if it is
absolutely continuous on each interval $[0,b]$ for $0 < b < +\infty$.
\end{definition}

\begin{remark}\label{rem-abs-cont}\rm\begin{enumerate} \item[(a)] It follows from the definition that an absolutely continuous function is differentiable almost 
everywhere, its derivative coincides with its distributional derivative almost everywhere and one can recover the function from its 
derivative $\dot x=y$ by the integration formula (i). 

\item[(b)] If $x:[0,b]\rightarrow {\cal H}$, where $b>0$, is absolutely continuous and $B:{\cal H}\rightarrow {\cal H}$ is 
$L$-Lipschitz continuous for $L\geq 0$, then the function $z=B\circ x$ is absolutely continuous, too.
This can be easily seen by using the characterization of absolute continuity in
Definition \ref{abs-cont}(iii). Moreover, $z$ is differentiable almost everywhere on $[0,b]$ and the inequality 
$\|\dot z (t)\|\leq L\|\dot x(t)\|$ holds for almost every $t \in [0,b]$.  
\end{enumerate}
\end{remark}

The following results, which can be interpreted as continuous counterparts of the quasi-Fej\'er monotonicity for sequences, will play 
an important role in the asymptotic analysis of the trajectories of the dynamical system investigated in this paper. 
For the proof of Lemma \ref{fejer-cont1} we refer the reader to \cite[Lemma 5.1]{abbas-att-sv}. 
Lemma \ref{fejer-cont2} follows by using similar arguments as used in  \cite[Lemma 5.2]{abbas-att-sv}.

\begin{lemma}\label{fejer-cont1} Suppose that $F:[0,+\infty)\rightarrow\R$ is locally absolutely continuous 
and bounded from below and that
there exists $G\in L^1([0,+\infty))$ such that for almost every $t \in [0,+\infty)$
\begin{equation*}\frac{d}{dt}F(t)\leq G(t).\end{equation*}
Then there exists $\lim_{t\rightarrow +\infty} F(t)\in\R$. 
\end{lemma}

\begin{lemma}\label{fejer-cont2}  If $1 \leq p < \infty$, $1 \leq r \leq \infty$, $F:[0,+\infty)\rightarrow[0,+\infty)$ is 
locally absolutely continuous, $F\in L^p([0,+\infty))\cap L^{\infty}([0,+\infty))$, $G_1,G_2:[0,+\infty)\rightarrow\R$, 
$G_1\in L^1([0,+\infty))$, $G_2\in  L^r([0,+\infty))$ and 
for almost every $t \in [0,+\infty)$
\begin{equation}\label{ineq-fejer-cont2}\frac{d}{dt}F(t)\leq G_1(t)+G_2(t),\end{equation} then $\lim_{t\rightarrow +\infty} F(t)=0$. 
\end{lemma}

\begin{proof} In case $r=1$ this follows from Lemma \ref{fejer-cont1} and the fact that $F\in L^p([0,+\infty))$. 

Assume now that $r>1$ and define $q:=1+p\left(1-\frac{1}{r}\right) >1$, which fulfills the relation 
\begin{equation}\label{conj-exp}\frac{q-1}{p}+\frac{1}{r}=1.
\end{equation}

Further, from \eqref{ineq-fejer-cont2} we derive for almost every $t \in [0,+\infty)$
\begin{equation}\label{ineq2-fejer-cont2}\frac{d}{dt}(F(t))^q\leq q(F(t))^{q-1}G_1(t)+q(F(t))^{q-1}G_2(t).\end{equation}
Since $F\in L^{\infty}([0,+\infty))$ and $G_1\in L^1([0,+\infty))$, the function $t\rightarrow (F(t))^{q-1}G_1(t)$ is 
$L^1$-integrable on $[0,+\infty)$. Moreover, due to $F^{q-1}\in L^\frac{p}{q-1}([0,+\infty))$, $G_2\in L^r([0,+\infty))$ and \eqref{conj-exp}, 
the function $t\rightarrow (F(t))^{q-1}G_2(t)$ is also $L^1$-integrable on $[0,+\infty)$. We conclude that 
the function on the right-hand side of inequality \eqref{ineq2-fejer-cont2} belongs to $L^1([0,+\infty))$. Applying now 
Lemma \ref{fejer-cont1} we obtain that there exists $\lim_{t\rightarrow +\infty} F(t)\in\R$, which combined again with 
$F\in L^p([0,+\infty))$ delivers the conclusion.  
\end{proof}

The next result which we recall here is the continuous version of the Opial Lemma. 

\begin{lemma}\label{opial-cont} Let $S \subseteq {\cal H}$ be a nonempty set and $x:[0,+\infty)\rightarrow{\cal H}$ a given function. Assume that 

(i) for every $x^*\in S$, $\lim_{t\rightarrow+\infty}\|x(t)-x^*\|$ exists; 

(ii) every weak sequential cluster point of the map $x$ belongs to $S$. 

\noindent Then there exists $x_{\infty}\in S$ such that  $x(t)$ converges weakly to $x_{\infty}$ as $t \rightarrow +\infty$. 
\end{lemma}

\section{A Levenberg-Marquardt dynamical system: existence and uniqueness of the trajectories}

Consider the optimization problem  
\begin{equation}\label{opt}\inf_{x\in\argmin\Psi}\{\Phi(x)+\Theta(x)\},\end{equation}
where ${\cal H}$ is a real Hilbert space and the following conditions hold:
\begin{align*}(H_\Psi)& \ \Psi:{\cal H}\rightarrow [0,+\infty) \mbox{ is convex, (Fr\'{e}chet) differentiable}
 \mbox{ with Lipschitz continuous gradient}\\& \mbox{ and }\argmin\Psi=\Psi^{-1}(0)\neq\emptyset;\\
 (H_\Theta)& \ \Theta:{\cal H}\rightarrow \R \mbox{ is convex, (Fr\'{e}chet) differentiable}
 \mbox{ with Lipschitz continuous gradient};\\
 (H_\Phi)& \ \Phi:{\cal H}\rightarrow \R\cup\{+\infty\} \mbox{ is convex, lower semicontinuous and fulfills the relation }\\
& \mbox{ }S:=\{z\in \argmin\Psi\cap\dom\Phi: \Phi(z)+\Theta(z)\leq \Phi(x)+\Theta(x) \  \forall x\in\argmin\Psi\}\neq\emptyset.
\end{align*}
Here, $\dom \Phi = \{x \in {\cal H}: \Phi(x) < +\infty\}$ denotes the effective domain of the function $\Phi$.

In connection with \eqref{opt}, we investigate the nonautonomous dynamical system  
\begin{equation}\label{dyn-syst}\left\{
\begin{array}{ll}
v(t)\in\partial\Phi(x(t))\\
\lambda(t)\dot x(t) + \dot v(t) + v(t) + \nabla \Theta(x(t))+\beta(t)\nabla \Psi(x(t))=0\\
x(0)=x_0, v(0)=v_0\in\partial \Phi(x_0),
\end{array}\right.\end{equation}
where $x_0,v_0\in {\cal H}$ and 
$$\partial \Phi : {\cal H} \rightrightarrows {\cal H}, \ \partial \Phi(x):=\{p\in {\cal H}: \Phi(y)\geq \Phi(x)+\langle p,y-x\rangle \ \forall y\in {\cal H}\},$$ for $\Phi(x) \in \R$ 
and $\partial \Phi(x):= \emptyset$ for $\Phi(x) \not\in \R$, denotes the convex subdifferential of $\Phi$. We denote by $\dom \partial \Phi = \{x \in {\cal H}: \partial \Phi(x) \neq \emptyset\}$ the domain of the operator $\partial \Phi$.

Furthermore, we make the following assumptions regarding the functions of time controlling the velocity and the penalty:
\begin{align*}(H^1_\lambda)& \ \lambda:[0,+\infty)\rightarrow(0,+\infty) \mbox{ is locally absolutely continuous};\\
 (H^1_\beta)& \ \beta:[0,+\infty)\rightarrow[0,+\infty)\mbox{ is locally integrable}. 
\end{align*}

Let us mention that due to $(H^1_\lambda)$, $\dot\lambda(t)$ exists for almost every $t\geq 0$.

\begin{remark}\rm\label{phi,psi0} (a) In case $\Phi(x)=0$ for all $x\in {\cal H}$, the dynamical system \eqref{dyn-syst} becomes 
\begin{equation}\label{dyn-syst-Phi0}\left\{
\begin{array}{ll}
\lambda(t)\dot x(t)  + \nabla \Theta(x(t))+\beta(t)\nabla \Psi(x(t))=0\\
x(0)=x_0,
\end{array}\right.\end{equation}
The asymptotic convergence of the trajectories generated by \eqref{dyn-syst-Phi0} has been investigated in \cite{att-cza-10} under the assumption $\lambda(t)=1$ for all $t\geq 0$, for $\Theta$ and $\Psi$ nonsmooth functions, by replacing their gradients with convex subdifferentials and, consequently, by treating the differential equation as a monotone inclusion (see \eqref{1ord-att-cz}).

(b) In case $\Psi(x)=0$ for all $x\in {\cal H}$, the dynamical system  
\begin{equation}\label{dyn-syst-Psi0}\left\{
\begin{array}{ll}
v(t)\in\partial\Phi(x(t))\\
\lambda(t)\dot x(t) + \dot v(t) + v(t) + \nabla \Theta(x(t))=0\\
x(0)=x_0, v(0)=v_0\in\partial \Phi(x_0),
\end{array}\right.\end{equation}
has been investigated in \cite{abbas-att-sv} (see, also, \cite{att-sv2011}, for the situation when $\Theta(x)=0$ for all $x\in {\cal H}$).

(c) In case $\Theta(x)=0$ and $\Psi(x)=\frac{1}{2}\|x\|^2$ for all $x\in {\cal H}$ and $\lambda(t) = 
\lambda \in \R$ for every $t\in[0,+\infty)$, the Levenberg-Marquardt dynamical system \eqref{dyn-syst} becomes 
\begin{equation}\label{dyn-syst-abbas}\left\{
\begin{array}{ll}
v(t)\in\partial\Phi(x(t))\\
\lambda\dot x(t) + \dot v(t) + v(t) +\beta(t)x(t)=0\\
x(0)=x_0, v(0)=v_0\in\partial \Phi(x_0).
\end{array}\right.\end{equation}
The dynamical system \eqref{dyn-syst-abbas} has been considered in \cite{abbas} 
in connection with the problem of finding the minimal norm elements among the minima of $\Phi$, namely, (see also \cite{att-cza-02} and \cite[Section 3]{att-maing})
\begin{equation}\label{opt2}\inf_{x\in\argmin\Phi}\|x\|^2.\end{equation} 
In contrast to \eqref{dyn-syst}, where the function describing the constrained set of \eqref{opt} is penalized, in \eqref{dyn-syst-abbas} 
the objective function of \eqref{opt2} is penalized via a vanishing penalization function (see \cite{abbas}). 
\end{remark}

In the following we specify what we understand under a solution of the dynamical system \eqref{dyn-syst}.

\begin{definition}\label{str-sol}\rm We say that the pair $(x,v)$ is a 
strong global solution of \eqref{dyn-syst}, if the following properties are satisfied: 
\begin{enumerate}
\item[(i)] $x,v:[0,+\infty)\rightarrow {\cal H}$ are locally absolutely continuous functions;

\item[(ii)] $v(t)\in\partial\Phi(x(t))$ for every $t\in[0,+\infty)$;

\item[(iii)] $\lambda(t)\dot x(t) + \dot v(t) + v(t) + \nabla \Theta(x(t))+\beta(t)\nabla \Psi(x(t))=0$ for almost every $t\in[0,+\infty)$;

\item[(iv)] $x(0)=x_0, v(0)=v_0$.
\end{enumerate}
\end{definition}

Similarly to the techniques used in \cite{att-sv2011}, we will show the existence and uniqueness of the trajectories generated by 
\eqref{dyn-syst} by converting it to an equivalent first order differential equation with respect to $z(\cdot)$, defined 
by 
\begin{equation}\label{def-z}z(t)=x(t)+\mu(t) v(t),\end{equation}
where $$\mu(t)=\frac{1}{\lambda(t)} \ \forall t\geq 0.$$ To this end we will make use of the resolvent and Yosida approximation 
of the convex subdifferential of $\Psi$. For $\gamma>0$, we denote by $$J_{\gamma\partial\Phi}=(\id+\gamma\partial\Phi)^{-1}$$ the resolvent 
of $\gamma\partial\Phi$. Due to the maximal monotonicity of $\partial\Phi$, the resolvent 
$J_{\gamma\partial\Phi}:{\cal H}\rightarrow{\cal H}$ is a 
single-valued operator with full-domain, which is, furthermore, nonexpansive, that is $1$-Lipschitz continuous. The Yosida regularization 
of $\partial\Phi$ is defined by 
$$(\partial\Phi)_{\gamma}=\frac{1}{\gamma}(\id-J_{\gamma\partial\Phi})$$
and it is $\gamma^{-1}$-Lipschitz continuous. For more properties of these operators we refer the reader to \cite{bauschke-book}. 

Assume now that $(x,v)$ is a strong global solution of \eqref{dyn-syst}. 
From \eqref{def-z} we have for every $t \in [0,+\infty)$
$$v(t)\in\partial\Phi(x(t))\Leftrightarrow z(t)\in x(t)+\mu(t)\partial\Phi(x(t))=(\id+\mu(t)\partial\Phi)(x(t)),$$
thus, from the definition of the resolvent we derive that relation (ii) in Definition \ref{str-sol} is equivalent to 
\begin{equation}\label{res-phi} x(t)=J_{\mu(t)\partial\Phi}(z(t)).\end{equation}
From \eqref{def-z}, \eqref{res-phi} and the definition of the Yosida regularization we obtain 
\begin{equation}\label{yos-phi}v(t)=(\partial\Phi)_{\mu(t)}(z(t)).\end{equation}

Further, by differentiating \eqref{def-z} and taking into account (iii) in Definition \ref{str-sol}, we get  for almost every $t\in[0,+\infty)$
\begin{align}\dot z(t)=\dot x(t)+\dot\mu(t)v(t)+\mu(t)\dot v(t) = \dot\mu(t)v(t)-\mu(t)v(t)-\mu(t)\nabla\Theta(x(t))-\beta(t)\mu(t)\nabla\Psi(x(t))\label{dot-z1}. 
\end{align}

Taking into account \eqref{res-phi}, \eqref{yos-phi} and \eqref{dot-z1} we conclude that $z$ defined in 
\eqref{def-z} is a strong global solution of the dynamical system  
\begin{equation}\label{dyn-syst-z}\left\{
\begin{array}{ll}
\!\!\dot z(t)+\big(\mu(t)-\dot\mu(t)\big)(\partial\Phi)_{\mu(t)}(z(t))+\mu(t)\nabla\Theta\Big(J_{\mu(t)\partial\Phi}(z(t))\Big)
+\beta(t)\mu(t)\nabla\Psi\Big(J_{\mu(t)\partial\Phi}(z(t))\Big)\!=\!0\\
z(0)=x_0+\mu(0)v_0.
\end{array}\right.\end{equation}

Vice versa, if $z$ is a strong global solution of \eqref{dyn-syst-z}, then one obtains via \eqref{res-phi} and \eqref{yos-phi} a strong global solution of \eqref{dyn-syst}.

\begin{remark}\rm\label{discrete} By considering the time discretization $\dot z(t)\approx\frac{z_{n+1}-z_n}{h_n}$ of the above dynamical system and by taking $\mu$ constant, from \eqref{res-phi} and 
\eqref{dyn-syst-z} we obtain the  iterative scheme 
\begin{equation}\label{discrete-syst-z}(\forall n \geq 0) \ \left\{
\begin{array}{ll}
x_n=J_{\mu\partial\Phi}(z_n)\\
z_{n+1}=(1-h_n)z_n+h_n\big(x_n-\mu\nabla\Theta(x_n)-\mu\beta_n\nabla\Psi(x_n)\big),
\end{array}\right.\end{equation}
which for $h_n=1$ yields the following algorithm 
\begin{equation}\label{discrete-syst-z2} (\forall n \geq 0) \ x_{n+1}= J_{\mu\partial\Phi}\big(x_n-\mu\nabla\Theta(x_n)-\mu\beta_n\nabla\Psi(x_n)\big).
\end{equation} 

The convergence of the above algorithm has been investigated in \cite{b-c-penalty-svva} in the more general 
framework of monotone inclusion problems, under the use of variable step sizes $(\mu_n)_{n\geq 0}$ and by assuming that
$$\forall p\in\ran N_{\argmin\Psi} \ \sum_{n \in \N} \mu_n\beta_n
\left[\Psi^*\left(\frac{p}{\beta_n}\right)-\sigma_{\argmin\Psi}\left(\frac{p}{\beta_n}\right)\right]<+\infty,$$
which is a condition that can be seen as a discretized version of the one stated in \eqref{ass-att-cz}. The case $\Theta(x)=0$ for all $x\in {\cal H}$ has been treated in \cite{att-cza-peyp-c} (see also the references therein). 
\end{remark}

Next we show that, given $x_0,v_0\in {\cal H}$ and by assuming $(H^1_{\lambda})$ and $(H^1_{\beta})$, there exists a unique strong global solution of the dynamical system \eqref{dyn-syst-z}. This will be done in the framework of the Cauchy-Lipschitz Theorem for absolutely continuous trajectories (see for example 
\cite[Proposition 6.2.1]{haraux}, \cite[Theorem 54]{sontag}). To this end we will make use of the following Lipschitz property 
of the resolvent operator as a function of the step size, which actually is a consequence of 
the classical results \cite[Proposition 2.6]{brezis} and \cite[Proposition 23.28]{bauschke-book} (see also \cite[Proposition 2.3]{att-sv2011} and \cite[Proposition 3.1]{abbas-att-sv}). 

\begin{proposition}\label{prop-Lipsch-res-step-size} Assume that $(H_\Phi)$ holds, $x\in {\cal H}$ and $0<\delta<+\infty$. 
Then the mapping $\tau\mapsto J_{\tau \partial\Phi}x$ is Lipschitz continuous on $[\delta,+\infty)$. More precisely, for any 
$\lambda_1,\lambda_2\in[\delta,+\infty)$ the following inequality holds:
\begin{equation}\label{Lipsch-res-step-size} \|J_{\lambda_1 \partial\Phi}x-J_{\lambda_2 \partial\Phi}x\|
\leq |\lambda_1-\lambda_2|\cdot\|(\partial\Phi)_\delta x\|.
\end{equation}
Furthermore, the function $\lambda\mapsto\|(\partial\Phi)_\lambda x\|$ is nonincreasing. 
\end{proposition}

Notice that the dynamical system \eqref{dyn-syst-z} can be written as 
\begin{equation}\label{existence}\left\{
\begin{array}{ll}
\dot z(t)=f(t,z(t))\\
z(0)=z_0,
\end{array}\right.\end{equation}
where $z_0=x_0+\mu(0)v_0$ and $f:[0,+\infty)\times {\cal H}\rightarrow {\cal H}$ is defined by 
\begin{equation}\label{f-existence}f(t,w)=
\big(\dot\mu(t)-\mu(t)\big)(\partial\Phi)_{\mu(t)}(w)-\mu(t)\nabla\Theta\Big(J_{\mu(t)\partial\Phi}(w)\Big)
-\beta(t)\mu(t)\nabla\Psi\Big(J_{\mu(t)\partial\Phi}(w)\Big).\end{equation}

In the following we denote by $L_{\nabla\Phi}$ and $L_{\nabla\Psi}$ the Lipschitz constants of $\nabla\Phi$ and $\nabla\Psi$, respectively. 

(a) Notice that for every $t\geq 0$ and every $w_1,w_2\in {\cal H}$ we have  

\begin{equation}\label{lipsch-x}\|f(t,w_1)-f(t,w_2)\|\leq 
\left(1+\frac{|\dot\lambda(t)|}{\lambda(t)}+\frac{L_{\nabla\Theta}}{\lambda(t)}+L_{\nabla\Psi}\frac{\beta(t)}{\lambda(t)}\right)\|w_1-w_2\|.\end{equation}
Indeed, this follows \eqref{f-existence}, the Lipschitz properties of the operators involved and the definition of $\mu(t)$. 
Further, notice that due to $(H^1_{\lambda})$ and $(H^1_{\beta})$, 
$$L_f: [0,+\infty) \rightarrow \R, 
L_f(t) = 1+\frac{|\dot\lambda(t)|}{\lambda(t)}+\frac{L_{\nabla\Theta}}{\lambda(t)}+L_{\nabla\Psi}\frac{\beta(t)}{\lambda(t)},$$
which is for every $t \geq 0$ equal to the Lipschitz-constant of $f(t,\cdot)$, satisfies
$$L_f(\cdot)\in L^1([0,b]) \mbox{ for any } 0<b<+\infty.$$ 

(b) We show now that 
\begin{equation}\label{existence-b}\forall w\in{\cal H}, \ \forall b>0, \ \ f(\cdot,w)\in L^1([0,b],{\cal H}).\end{equation}

We fix $w\in{\cal H}$ and $b>0$. Due to $(H^1_{\lambda})$, there exist $\lambda_{\min}, \lambda_{\max}>0$ such that
$$0<\lambda_{\min}\leq\lambda(t)\leq  \lambda_{\max}\ \forall t\in[0,b],$$
hence $$0<\frac{1}{\lambda_{\max}}\leq\mu(t)\leq  \frac{1}{\lambda_{\min}}\ \forall t\in[0,b].$$

Relying on Proposition \ref{prop-Lipsch-res-step-size} we obtain for all $t\in[0,b]$ the following chain of inequalities: 
\begin{align*}
\|f(t,w)\| \leq & \ |\dot\mu(t)-\mu(t)|\cdot\|(\partial\Phi)_{\frac{1}{\lambda_{\max}}}(w)\|\\
 & \ +\mu(t)\|\nabla\Theta(J_{\mu(t)\partial\Phi}(w))-\nabla\Theta(J_{\frac{1}{\lambda_{\max}}\partial\Phi}(w))\| +\mu(t)\|\nabla\Theta(J_{\frac{1}{\lambda_{\max}}\partial\Phi}(w))\|\\
& \ +\beta(t)\mu(t)\|\nabla\Psi(J_{\mu(t)\partial\Phi}(w))-\nabla\Psi(J_{\frac{1}{\lambda_{\max}}\partial\Phi}(w))\| +\beta(t)\mu(t)\|\nabla\Psi(J_{\frac{1}{\lambda_{\max}}\partial\Phi}(w))\|\\
\leq & \ |\dot\mu(t)-\mu(t)|\cdot\|(\partial\Phi)_{\frac{1}{\lambda_{\max}}}(w)\|\\
 & \ +L_{\nabla\Theta}\mu(t)\left(\mu(t)-\frac{1}{\lambda_{\max}}\right)\cdot\|\nabla\Theta(J_{\frac{1}{\lambda_{\max}}\partial\Phi}(w))\|+\mu(t)\|\nabla\Theta(J_{\frac{1}{\lambda_{\max}}\partial\Phi}(w))\|\\
& \ +L_{\nabla\Psi}\beta(t)\mu(t)\left(\mu(t)-\frac{1}{\lambda_{\max}}\right)\cdot\|\nabla\Psi(J_{\frac{1}{\lambda_{\max}}\partial\Phi}(w))\| +\beta(t)\mu(t)\|\nabla\Psi(J_{\frac{1}{\lambda_{\max}}\partial\Phi}(w))\|.
\end{align*}
Now \eqref{existence-b} follows from the properties of the functions $\mu$ and $\beta$, and the fact that 
$$|\dot\mu(t)-\mu(t)|\leq \frac{1}{\lambda_{\min}}\left(1+\frac{|\dot\lambda(t)|}{\lambda_{\min}}\right) \mbox{ for almost every }t\geq 0.$$

In the light of the statements proven in (a) and (b), the existence and uniqueness of a strong global solution of the dynamical system 
\eqref{dyn-syst-z} follow from \cite[Proposition 6.2.1]{haraux} (see also \cite[Theorem 54]{sontag}).

Finally, similarly to the proof of \cite[Theorem 2.4(ii)]{att-sv2011}, one can guarantee the existence and uniqueness of the trajectories 
generated by \eqref{dyn-syst} by relying on the properties of the dynamical system \eqref{dyn-syst-z} and on \eqref{res-phi} and \eqref{yos-phi}. The details are left to the reader.

\section{Convergence of the trajectories and of the objective function values}\label{sec3}

In this section we prove weak convergence for the trajectory generated by the dynamical system \eqref{dyn-syst} to an 
optimal solution of \eqref{opt} as well as convergence for the objective function values of the latter along the trajectory. Some techniques 
from \cite{att-cza-10} and \cite{att-sv2011} will be useful in this context. 

To this end we will make the following supplementary assumptions:
\begin{align*}(H^2_\lambda)& \ \lambda:[0,+\infty)\rightarrow(0,+\infty) 
\mbox{ is locally absolutely continuous and }\dot\lambda(t)\leq 0\mbox{ for almost every }\\
& \ t \in [0,+\infty);\\
 (H^2_\beta)& \ \beta:[0,+\infty)\rightarrow(0,+\infty)\mbox{ is measurable and bounded from above on each interval }[0,b],\\
& \ 0<b<+\infty;\\ 
 (H)& \ \forall p\in\ran N_{\argmin \Psi} \ \int_0^{+\infty} 
\beta(t)\left[\Psi^*\left(\frac{p}{\beta(t)}\right)-\sigma_{\argmin \Psi}\left(\frac{p}{\beta(t)}\right)\right]dt<+\infty; \\
(\widetilde{H})& \ \partial(\Phi+\Theta+\delta_{\argmin\Psi})=\partial\Phi+\nabla\Theta+N_{\argmin\Psi}, 
\end{align*}
where \begin{itemize}
\item $N_{\argmin\Psi}$ is the normal cone to the set $\argmin\Psi$: $N_{\argmin\Psi}(x)=\{p\in{\cal H}:\langle p,y-x\rangle\leq 0 \ 
\forall y\in \argmin\Psi\}$ for  $x \in \argmin\Psi$ and $N_{\argmin\Psi}(x)=\emptyset$ for $x\not\in \argmin\Psi$; 
\item $\ran N_{\argmin\Psi}$ is the range of the normal cone $N_{\argmin\Psi}$: $p\in\ran N_{\argmin\psi}$ if and only if 
there exists $x\in \argmin\Psi$ such that $p\in N_{\argmin\Psi}(x)$; 
\item $\Psi^*: {\cal H}\rightarrow \R\cup\{+\infty\}$ is the Fenchel conjugate of $\Psi$:
$\Psi^*(p)=\sup_{x\in{\cal H}}\{\langle p,x\rangle-\Psi(x)\} \ \forall p\in {\cal H};$ 
\item $\sigma_{\argmin\Psi}: {\cal H}\rightarrow \R\cup\{+\infty\}$ is the support function of the set $\argmin\Psi$:  
$\sigma_{\argmin\Psi}(p)=\sup_{x\in {\argmin\Psi}}\langle p,x\rangle$ for all $p\in {\cal H}$; 
\item $\delta_{\argmin\Psi}:{\cal H}\rightarrow\R\cup\{+\infty\}$ is 
the indicator function of ${\argmin\Psi}$: it takes the value $0$ on the set ${\argmin\Psi}$ and $+\infty$, otherwise.
\end{itemize}

We have $N_{\argmin\Psi}=\partial \delta_{\argmin\Psi}$. Moreover, $p\in N_{\argmin\Psi}(x)$ if and only if  $x\in {\argmin\Psi}$ and $\sigma_{\argmin\Psi}(p)=\langle p,x\rangle$. 

\begin{remark}
\begin{itemize} 
\item[(a)] The condition $\dot\lambda(t)\leq 0$ for almost every $t \in [0,+\infty)$ has been used in \cite{att-sv2011} in the study of the asymptotic convergence of the dynamical system \eqref{att-sv-dyn}, 
when approaching the monotone inclusion problem \eqref{o-tx}. 
\item[(b)] Under $(H_\Psi)$, due to $\Psi\leq\delta_{\argmin \Psi}$, we 
have $\Psi^*\geq\delta_{\argmin \Psi}^*=\sigma_{\argmin \Psi}.$
\item[(c)] When $\Psi=0$ (see Remark \ref{phi,psi0}(b)), it holds $N_{\argmin \Psi}(x)=\{0\}$ for every $x\in \argmin \Psi={\cal H}$, 
$\Psi^*=\sigma_{\argmin \Psi}=\delta_{\{0\}}$, which shows that in this case $(H)$ trivially holds. 
\item[(d)] A nontrivial situation in which condition $(H)$ is fulfilled is when $\psi(x)=\frac{1}{2}\inf_{y\in C}\|x-y\|^2$, for a nonempty, convex and closed set $C \subseteq {\cal H}$ (see \cite{att-cza-10}). Then \eqref{ass-att-cz} holds if and only if 
$$\int_0^{+\infty}\frac{1}{\beta(t)} dt<+\infty,$$ which is trivially satisfied for $\beta(t)=(1+t)^\alpha$ with $\alpha>1$. 
\item[(e)] Due to the continuity of $\Theta$, the condition $(\widetilde{H})$ is equivalent to  
$$\partial(\Phi+\delta_{\argmin\Psi})=\partial\Phi+N_{\argmin\Psi},$$ which holds when
$0\in \sqri(\dom\Phi-\argmin\Psi)$, a condition that is fulfilled, if $\Phi$ is continuous at a point in 
$\dom \Phi\cap \argmin\Psi$ or $\inte (\argmin\Psi)\cap\dom \Phi\neq\emptyset$ (we invite the reader to consult also \cite{bauschke-book}, \cite{b-hab} and \cite{Zal-carte} for 
other sufficient conditions for the above subdifferential sum formula). Here, for $M\subseteq {\cal H}$ a convex set, 
 $$\sqri M:=\{x\in M:\cup_{\lambda>0}\lambda(M-x) \ \mbox{is a closed linear subspace of} \ {\cal H}\}$$
denotes its strong quasi-relative interior. We always have $\inte M\subseteq\sqri M$ (in general this inclusion may be strict). 
If ${\cal H}$ is finite-dimensional, then $\sqri M$ coincides with $\ri M$, the relative interior of $M$, which is the interior of $M$ with 
respect to its affine hull.
\end{itemize}
\end{remark}

The following differentiability result of the composition of convex functions with absolutely 
continuous trajectories that is due to Br\'{e}zis (see \cite[Lemme 4, p. 73]{brezis} and also \cite[Lemma 3.2]{att-cza-10}) will play an important role in our analysis. 

\begin{lemma}\label{diff-brezis} Let $f:{\cal H}\rightarrow \R\cup\{+\infty\}$ be a proper, convex and lower semicontinuous function. 
Let $x\in L^2([0,T],{\cal H})$ be absolutely continuous such that $\dot x\in L^2([0,T],{\cal H})$ and $x(t)\in\dom f$ for almost every 
$t \in [0,T]$. Assume that there exists $\xi\in L^2([0,T],{\cal H})$ such that $\xi(t)\in\partial f(x(t))$ for almost every 
$t \in [0,T]$. Then the function $t\mapsto f(x(t))$ is absolutely continuous and for every $t$ such that $x(t)\in\dom \partial f$ 
we have $$\frac{d}{dt}f(x(t))=\langle \dot x(t),h\rangle \ \forall h\in\partial f(x(t)).$$ 
\end{lemma}

We start our convergence analysis with the following technical result.

\begin{lemma}\label{diff-Phix} Assume that $(H_\Psi)$, $(H_\Theta)$, $(H_\Phi)$, $(H^1_{\lambda})$ and $(H^2_{\beta})$ hold 
and let $(x,v):[0,+\infty)\rightarrow{\cal H} \times {\cal H}$ be  a strong stable solution of the dynamical system 
\eqref{dyn-syst}. Then the following statements are true: 
\begin{enumerate}
 \item[(i)] $\langle \dot x(t),\dot v(t)\rangle\geq 0$ for almost every $t\in [0,+\infty)$;  
 \item[(ii)] $\frac{d}{dt}\Phi(x(t))=\langle \dot x(t),v(t)\rangle$ for almost every $t\in [0,+\infty)$.
 \end{enumerate}
\end{lemma}

\begin{proof} (i) See \cite[Proposition 3.1]{att-sv2011}. The proof relies on the first relation in \eqref{dyn-syst} and 
the monotonicity of the convex subdifferential. 

(ii) The proof makes use of Lemma \ref{diff-brezis}. Let $T>0$ be fixed. Due to the continuity of $x$ and $v$ we obviously have 
$$x,v\in L^2([0,T],{\cal H}).$$ The only condition which has to be checked is $\dot x \in L^2([0,T],{\cal H})$. 
By considering the second relation in \eqref{dyn-syst} and by inner multiplying it with $\dot x(t)$, we derive for almost every $t\in[0,T]$
$$\lambda(t)\|\dot x(t)\|^2+\langle \dot x(t),\dot v(t)\rangle+\langle \dot x(t),v(t)\rangle+ \langle \dot x(t),\nabla\Theta(x(t))\rangle + \beta(t)\langle \dot x(t),\nabla\Psi(x(t))\rangle=0.$$
Using (i) we obtain for almost every $t\in[0,T]
$\begin{equation}\label{dx-l2-loc}\lambda(t)\|\dot x(t)\|^2+\langle \dot x(t),v(t)\rangle+\langle \dot x(t),\nabla\Theta(x(t))\rangle+
\beta(t)\langle \dot x(t),\nabla\Psi(x(t))\rangle\leq 0.\end{equation}
Since $x,v$ are continuous on $[0,T]$, they are bounded on $[0,T]$, a property which is shared also by 
$t\mapsto \beta(t)\nabla\Psi(x(t))$, due to $(H^2_{\beta})$ and $(H_\Psi)$, and by $t\mapsto \nabla\Theta(x(t))$, due to $(H_\Theta)$. Since $\lambda$ is bounded from below by a positive 
constant on $[0,T]$, from \eqref{dx-l2-loc} one easily obtains that $$\dot x \in L^2([0,T],{\cal H})$$
and the conclusion follows by applying Lemma \ref{diff-brezis}.
\end{proof}

\begin{lemma}\label{l1-op} Assume that $(H_\Psi)$, $(H_\Theta)$, $(H_\Phi)$, $(H^2_{\lambda})$, $(H^2_{\beta})$, $(H)$ 
and $(\widetilde{H})$ hold and let $(x,v):[0,+\infty)\rightarrow{\cal H} \times {\cal H}$ be a strong stable solution of the dynamical system 
\eqref{dyn-syst}. Choose arbitrary $z\in S$ and $p\in N_{\argmin\Psi}(z)$ such that $-p-\nabla\Theta(z)\in\partial\Phi(z)$. 
Define $g_z,h_z:[0,+\infty)\rightarrow [0,+\infty)$ as
$$g_z(t)=\Phi(z)-\Phi(x(t))+\langle v(t),x(t)-z\rangle$$ and 
$$h_z(t)=\Theta(z)-\Theta(x(t))+\langle \nabla\Theta(x(t)),x(t)-z\rangle.$$
The following statements are true: 
\begin{enumerate}
 \item[(i)] $\exists\lim_{t\rightarrow+\infty}\left(\frac{\lambda(t)}{2}\|x(t)-z\|^2+g_z(t)\right)\in[0,+\infty)$; 
 \item[(ii)] $\int_0^{+\infty}\beta(t)\Psi(x(t))dt<+\infty$; 
 \item[(iii)] $\exists\lim_{t\rightarrow+\infty}\int_0^t\langle p,x(s)-z\rangle ds\in\R$;
 \item[(iv)] $\exists\lim_{t\rightarrow+\infty}\int_0^t\Big((\Phi+\Theta)(x(s))-(\Phi+\Theta)(z)+\beta(s)\Psi(x(s))\Big)ds\in\R$; 
 \item[(v)] $\exists\lim_{t\rightarrow+\infty}\int_0^t\big(\langle v(s),x(s)-z\rangle +
 \langle \nabla\Theta(x(s)),x(s)-z\rangle+\beta(s)\Psi(x(s)\big)ds\in\R$;
 \item[(vi)] $\exists\lim_{t\rightarrow+\infty}\int_0^t\Big((\Phi+\Theta)(x(s))-(\Phi+\Theta)(z)\Big)ds\in\R$;
 \item[(vii)] $\exists\lim_{t\rightarrow+\infty}\int_0^t\big(\langle v(s),x(s)-z\rangle +
 \langle \nabla\Theta(x(s)),x(s)-z\rangle\big)ds\in\R$;
 \item[(viii)] $g_z\in L^1([0,+\infty))\cap L^\infty([0,+\infty))$ and $h_z\in L^1([0,+\infty))$. 
\end{enumerate}
\end{lemma}

\begin{proof} For the beginning, we notice that from the definition of $S$ and $(\widetilde{H})$ we have 
$$0\in\partial(\Phi+\Theta+\delta_{\argmin\Psi})(z)=\partial\Phi(z)+\nabla\Theta(z)+N_{\argmin\Psi}(z),$$
hence there exists \begin{equation}\label{p-in-n}p\in N_{\argmin\Psi}(z)\end{equation} such that 
\begin{equation}\label{-p-nabla-in}-p-\nabla\Theta(z)\in\partial\Phi(z).\end{equation}
For almost every $t\geq 0$ it holds according to \eqref{dyn-syst} 
\begin{align}\label{dt}\frac{d}{dt}\left(\frac{\lambda(t)}{2}\|x(t)-z\|^2\right)
& = \ \frac{\dot\lambda(t)}{2}\|x(t)-z\|^2 + \lambda(t)\langle \dot x(t),x(t)-z\rangle \nonumber\\
& = \  \frac{\dot\lambda(t)}{2}\|x(t)-z\|^2 - \langle \dot v(t),x(t)-z\rangle -\langle v(t),x(t)-z\rangle\nonumber\\
& \ \ \ \ \ -\langle \nabla\Theta(x(t)),x(t)-z\rangle-\beta(t)\langle  \nabla\Psi(x(t)),x(t)-z\rangle.
\end{align}
From \eqref{dyn-syst} and the convexity of $\Phi, \Theta$ and $\Psi$ we have for every $t \in [0,+\infty)$
\begin{equation}\label{Phi-conv} \Phi(z)\geq\Phi(x(t))+\langle v(t),z-x(t)\rangle
\end{equation}
\begin{equation}\label{Theta-conv} \Theta(z)\geq\Theta(x(t))+\langle \nabla\Theta(x(t)),z-x(t)\rangle
\end{equation}
and 
\begin{equation}\label{Psi-conv} 0=\Psi(z)\geq\Psi(x(t))+\langle \nabla\Psi(x(t)),z-x(t)\rangle.
\end{equation}
From \eqref{-p-nabla-in} and the convexity $\Phi$ and $\Theta$ we obtain for every $t \in [0,+\infty)$
\begin{equation}\label{Phi-conv2}\Phi(x(t))\geq\Phi(z)+\langle -p-\nabla\Theta(z),x(t)-z\rangle
\end{equation}
and 
\begin{equation}\label{Theta-conv2}\Theta(x(t))\geq\Theta(z)+\langle \nabla\Theta(z),x(t)-z\rangle.
\end{equation}

Further, due to Lemma \ref{diff-Phix}(ii) it holds for almost every $t \in [0,+\infty)$
\begin{align}\label{dt-gz}\frac{d}{dt}g_z(t)=& \ -\langle \dot x(t),v(t)\rangle+\langle \dot v(t),x(t)-z\rangle+
\langle v(t),\dot x(t)\rangle\nonumber\\ = & \ \langle \dot v(t),x(t)-z\rangle. 
\end{align}

On the other hand, using \eqref{p-in-n} and the Young-Fenchel inequality we obtain for every $t \in [0,+\infty)$
\begin{align}\label{ineq-conj}\beta(t)\Psi(x(t))+\langle -p,x(t)-z\rangle & =
\beta(t)\left(\Psi(x(t))+\left\langle \frac{-p}{\beta(t)},x(t)-z\right\rangle\right)\nonumber\\
& = \beta(t)\left(\Psi(x(t))-\left\langle \frac{p}{\beta(t)},x(t)\right\rangle+\sigma_{\argmin \Psi}\left(\frac{p}{\beta(t)}\right)\right)\nonumber\\
& \geq \beta(t)\left(-\Psi^*\left(\frac{p}{\beta(t)}\right)+\sigma_{\argmin \Psi}\left(\frac{p}{\beta(t)}\right)\right).\end{align}
Finally, we obtain for almost every $t \in [0,+\infty)$
\begin{align} &  \frac{d}{dt}\left(\frac{\lambda(t)}{2}\|x(t)-z\|^2+g_z(t)\right)
+\beta(t)\left(-\Psi^*\left(\frac{p}{\beta(t)}\right)+\sigma_{\argmin \Psi}\left(\frac{p}{\beta(t)}\right)\right)\nonumber\\
\leq & \ \frac{d}{dt}\left(\frac{\lambda(t)}{2}\|x(t)-z\|^2+g_z(t)\right)+\beta(t)\Psi(x(t))+\langle -p,x(t)-z\rangle\nonumber\\
\leq & \ \frac{d}{dt}\left(\frac{\lambda(t)}{2}\|x(t)-z\|^2+g_z(t)\right)
+(\Phi+\Theta)(x(t))-(\Phi+\Theta)(z)+\beta(t)\Psi(x(t))\nonumber\\
\leq & \ \frac{d}{dt}\left(\frac{\lambda(t)}{2}\|x(t)-z\|^2+g_z(t)\right)+
\langle v(t),x(t)-z\rangle+\langle \nabla\Theta(x(t)),x(t)-z\rangle+\beta(t)\Psi(x(t))\nonumber\\
\leq & \ 0,\label{ineq3}
\end{align}
where the first inequality follows from \eqref{ineq-conj}, the second one from \eqref{Phi-conv2} and \eqref{Theta-conv2}, 
the next one from \eqref{Phi-conv} and \eqref{Theta-conv}, and the last one from $(H^2_{\lambda})$, \eqref{dt}, \eqref{dt-gz} and \eqref{Psi-conv}.  

(i) Since for almost every $t \in [0,+\infty)$ we have (see \eqref{ineq3})
$$\frac{d}{dt}\left(\frac{\lambda(t)}{2}\|x(t)-z\|^2+g_z(t)\right)
\leq \beta(t)\left(\Psi^*\left(\frac{p}{\beta(t)}\right)-\sigma_{\argmin \Psi}\left(\frac{p}{\beta(t)}\right)\right),$$
the conclusion follows from Lemma \ref{fejer-cont1}, $(H)$ and the fact that $g_z(t)\geq 0$ for every $t\geq 0$. 

(ii) Let $F:[0,+\infty)\rightarrow\R$ be defined by 
$$F(t)=\int_0^t\big(-\beta(s)\Psi(x(s))+\langle p,x(s)-z\rangle\big)ds \ \forall t \in [0,+\infty).$$
From \eqref{ineq3} we have for almost every $s \in [0,+\infty)$
$$-\beta(s)\Psi(x(s))+\langle p,x(s)-z\rangle\geq \frac{d}{ds}\left(\frac{\lambda(s)}{2}\|x(s)-z\|^2+g_z(s)\right).$$
By integration we obtain for every $t \in [0,+\infty)$
\begin{align*}F(t)\geq & \ \frac{\lambda(t)}{2}\|x(t)-z\|^2+g_z(t)-\frac{\lambda(0)}{2}\|x_0-z\|^2-g_z(0)\\
\geq&  \  -\frac{\lambda(0)}{2}\|x_0-z\|^2-g_z(0), 
\end{align*}
hence $F$ is bounded from below. Furthermore, from \eqref{ineq-conj} we derive  for every $t \in [0,+\infty)$
\begin{align*}\frac{d}{dt}F(t)= & \ -\beta(t)\Psi(x(t))+\langle p,x(t)-z\rangle\\
\leq & \ \beta(t)\left(\Psi^*\left(\frac{p}{\beta(t)}\right)-\sigma_{\argmin \Psi}\left(\frac{p}{\beta(t)}\right)\right). 
\end{align*}
From $(H)$ and Lemma \ref{fejer-cont1} it follows that $\lim_{t\rightarrow+\infty}F(t)$ exists and it is a real number. 
Hence 
\begin{equation}\label{b}\exists \lim_{t\rightarrow+\infty}\int_0^t\big(\beta(s)\Psi(s)+
\langle-p,x(s)-z\rangle\big)ds\in\R.\end{equation}
Further, since $\psi\geq 0$, we obtain for every $t \in [0,+\infty)$ 
$$\beta(t)\Psi(x(t))+\langle -p,x(t)-z\rangle\geq \frac{\beta(t)}{2}\Psi(x(t))+\langle -p,x(t)-z\rangle.$$ 
Similarly to \eqref{ineq-conj}one can show that for every $t \in [0,+\infty)$ 
$$\frac{\beta(t)}{2}\Psi(x(t))+\langle -p,x(t)-z\rangle
\geq \frac{\beta(t)}{2}\left(-\Psi^*\left(\frac{2p}{\beta(t)}\right)+\sigma_{\argmin \Psi}\left(\frac{2p}{\beta(t)}\right)\right),$$
while from \eqref{ineq3} we obtain that for  almost every $t \in [0,+\infty)$ it holds
\begin{align*} &  \frac{d}{dt}\left(\frac{\lambda(t)}{2}\|x(t)-z\|^2+g_z(t)\right)
+\frac{\beta(t)}{2}\left(-\Psi^*\left(\frac{2p}{\beta(t)}\right)+\sigma_{\argmin \Psi}\left(\frac{2p}{\beta(t)}\right)\right)\nonumber\\
\leq & \ \frac{d}{dt}\left(\frac{\lambda(t)}{2}\|x(t)-z\|^2+g_z(t)\right)+\frac{\beta(t)}{2}\Psi(x(t))+\langle -p,x(t)-z\rangle\nonumber\\
\leq & \ \frac{d}{dt}\left(\frac{\lambda(t)}{2}\|x(t)-z\|^2+g_z(t)\right)
+\beta(t)\Psi(x(t))+\langle -p,x(t)-z\rangle\nonumber\\
\leq & \ 0
\end{align*}
By using the same arguments as used in the proof of \eqref{b} it yields that
\begin{equation}\label{b2}\exists \lim_{t\rightarrow+\infty}\int_0^t\left(\frac{\beta(s)}{2}\Psi(s)+\langle -p,x(s)-z\rangle\right)ds\in\R.\end{equation}
Finally, from \eqref{b} and \eqref{b2} we obtain (ii). 

(iii) Follows from \eqref{b} and (ii). 

(iv)-(v) These statements follow from \eqref{ineq3} and \eqref{ineq-conj}, by using similar arguments as used for 
proving \eqref{b}.  

(vi)-(vii) These statements are direct consequences of (iv), (v) and (ii). 

(viii) Combining (vi) and (vii) with $g_z,h_z\geq 0$, we easily derive that $$g_z+h_z\in L^1([0,+\infty)).$$
Since $$0\leq g_z\leq g_z+h_z,$$ we deduce that $g_z\in L^1([0,+\infty))$ and $h_z\in L^1([0,+\infty))$. 
Finally, notice that due to (i) there exists $T >0$ such that $g_z$ is bounded on $[T, +\infty)$. The boundedness of $g_z$ on
$[0,T]$ follows from \eqref{Phi-conv2} and the continuity of $x$ and $v$. Thus, $g_z\in L^\infty([0,+\infty))$.
\end{proof}

In order to proceed with the asymptotic analysis of the dynamical system \eqref{dyn-syst}, we make the following more involved assumptions on the functions $\lambda$ and $\beta$, respectively:
\begin{align*}(H^3_\lambda)& \ \lambda:[0,+\infty)\rightarrow(0,+\infty) 
\mbox{ is locally absolutely continuous}, \dot\lambda(t)\leq 0\mbox{ for almost every }\\
& \ t\in [0,+\infty)\mbox{ and }\lim_{t\rightarrow+\infty}\lambda(t)>0;\\
 (H^3_\beta)& \ \beta :[0,+\infty) \rightarrow (0,+\infty) \mbox{ is locally absolutely continuous, it satisfies for some } k\geq 0\mbox{ the}\\ 
 & \ \mbox{growth condition } 0\leq\dot\beta(t)\leq k\beta(t) \mbox{ for almost every }t \in [0, +\infty)\mbox{ and } \lim_{t\rightarrow+\infty}\beta(t)=+\infty.
 \end{align*}
 
\begin{lemma}\label{l-psix=0} Assume that $(H_\Psi)$, $(H_\Theta)$, $(H_\Phi)$, $(H^3_{\lambda})$, $(H^3_{\beta})$, $(H)$ 
and $(\widetilde{H})$ hold and let $(x,v):[0,+\infty)\rightarrow{\cal H} \times {\cal H}$ be a strong stable solution of the dynamical system 
\eqref{dyn-syst}. The following statements are true: 
\begin{enumerate}
 \item[(i)] $x$ is bounded; 
 \item[(ii)] $\lim_{t\rightarrow+\infty}\Psi(x(t))=0$.
\end{enumerate}
\end{lemma}

\begin{proof} Take an arbitrary $z\in S$ and (according to $(\widetilde H)$) $p\in N_{\argmin\Psi}(z)$ such that $-p-\nabla\Theta(z)\in\partial\Phi(z)$ and 
consider the functions $g_z,h_z$ defined in Lemma \ref{l1-op}. 

(i) According to Lemma \ref{l1-op}(i), since $g_z\geq 0$, we have 
that $t\mapsto\lambda(t)\|x(t)-z\|^2$ is bounded, which combined with $\lim_{t\rightarrow+\infty}\lambda(t)>0$ implies that 
$x$ is bounded. 

(ii) Consider the function $E_1:[0,+\infty)\rightarrow \R$ defined for every $t\in [0,+\infty)$ by 
$$E_1(t)=\frac{(\Phi+\Theta)(x(t))}{\beta(t)}+\Psi(x(t)).$$
Using Lemma \ref{diff-Phix} and \eqref{dyn-syst} we obtain for almost every $t\in [0,+\infty)$
\begin{align} \dot E_1(t)= & \ \frac{1}{\beta(t)}\big(\langle v(t),\dot x(t)\rangle+\langle \nabla\Theta(x(t)),\dot x(t)\rangle\big)
-\frac{\dot\beta(t)}{\beta^2(t)}(\Phi+\Theta)(x(t))+\langle \nabla\Psi(x(t)),\dot x(t)\rangle\nonumber\\
= & \ \frac{1}{\beta(t)}\langle v(t)+\nabla\Theta(x(t))+\beta(t)\nabla\Psi(x(t)),\dot x(t)\rangle
-\frac{\dot\beta(t)}{\beta^2(t)}(\Phi+\Theta)(x(t))\nonumber\\
= & \ \frac{1}{\beta(t)}\langle -\lambda(t)\dot x(t)-\dot v(t),\dot x(t)\rangle
-\frac{\dot\beta(t)}{\beta^2(t)}(\Phi+\Theta)(x(t))\nonumber\\
= & \ -\frac{\lambda(t)}{\beta(t)}\|\dot x(t)\|^2 -\frac{1}{\beta(t)}\langle\dot v(t),\dot x(t)\rangle
-\frac{\dot\beta(t)}{\beta^2(t)}(\Phi+\Theta)(x(t))\nonumber\\
\leq & \ -\frac{\dot\beta(t)}{\beta^2(t)}\inf_{t\geq 0}(\Phi+\Theta)(x(t))\label{ineq-E1}, 
\end{align}
where we used that, according to \eqref{Phi-conv2}, \eqref{Theta-conv2} and (i), $(\Phi+\Theta)(x(t))$ is bounded from below. 
From \eqref{ineq-E1} and Lemma \ref{fejer-cont1} it follows that there exists $\lim_{t\rightarrow+\infty}E_1(t)\in\R$.

Using now Lemma \ref{l1-op}(iv) we get 
\begin{equation}\label{liminf}\liminf_{t\rightarrow+\infty}\Big((\Phi+\Theta)(x(t))-(\Phi+\Theta)(z)+
\beta(t)\Psi(x(t))\Big)\leq 0\end{equation}
and, since $(\Phi+\Theta)(x(t))$ is bounded from below, this limes inferior  is a real number. Let $(t_n)_{n\in\N}$ be a 
sequence with $\lim_{n\rightarrow+\infty}t_n=+\infty$ such that 
\begin{align*}
& \lim_{n\rightarrow+\infty} \Big((\Phi+\Theta)(x(t_n))-(\Phi+\Theta)(z)+
\beta(t_n)\Psi(x(t_n))\Big) =\\
& \liminf_{t\rightarrow+\infty}\Big((\Phi+\Theta)(x(t))-(\Phi+\Theta)(z)+\beta(t)\Psi(x(t))\Big)\in\R.
\end{align*}

Since
$$E_1(t_n)=\frac{1}{\beta(t_n)}\Big((\Phi+\Theta)(x(t_n))-(\Phi+\Theta)(z)+\beta(t_n)\Psi(x(t_n))\Big)
+\frac{(\Phi+\Theta)(z)}{\beta(t_n)}  \ \forall n \in \N$$
and $\lim_{n\rightarrow+\infty}\beta(t_n)=+\infty$, it yields that $\lim_{n\rightarrow+\infty}E_1(t_n)=0$. 
Thus, since $\lim_{t\rightarrow+\infty}E_1(t)$ exists,
$$\lim_{t\rightarrow+\infty}E_1(t)=0.$$

The statement follows by taking into consideration that for every $t \in [0,+\infty)$
$$0\leq\Psi(x(t))\leq \Psi(x(t))+\frac{1}{\beta(t)}\Big((\Phi+\Theta)(x(t))-\inf_{s\geq 0}(\Phi+\Theta)(x(s))\Big)=
E_1(t)-\frac{1}{\beta(t)}\inf_{s\geq 0}(\Phi+\Theta)(x(s))$$
in combination with $\lim_{t\rightarrow+\infty}\beta(t)=+\infty$.
\end{proof}

\begin{lemma}\label{l-att-cz} Assume that $(H_\Psi)$, $(H_\Theta)$, $(H_\Phi)$, $(H^3_{\lambda})$, $(H^3_{\beta})$, $(H)$ 
and $(\widetilde{H})$ hold and let $(x,v):[0,+\infty)\rightarrow{\cal H} \times {\cal H}$ be a strong stable solution of the dynamical system 
\eqref{dyn-syst}. Then
$$\liminf_{t\rightarrow+\infty}(\Phi+\Theta)(x(t))\geq (\Phi+\Theta)(z) \ \forall z\in S.$$
\end{lemma}

\begin{proof} Take an arbitrary $z\in S$. From $(\widetilde{H})$ there exists $p\in N_{\argmin\Psi}(z)$ such that 
$-p-\nabla\Theta(z)\in\partial\Phi(z)$. From Lemma \ref{l1-op}(iii) we get 
\begin{equation}\label{liminf-p}\liminf_{t\rightarrow+\infty}\langle -p,x(t)-z\rangle\leq 0.\end{equation} We claim that 
\begin{equation}\label{liminf-p-0} \liminf_{t\rightarrow+\infty}\langle -p,x(t)-z\rangle= 0.\end{equation}
Since according to the previous lemma $x$ is bounded, this limit inferior is a real number. Let $(t_n)_{n\in\N}$ be a 
sequence with $\lim_{n\rightarrow+\infty}t_n=+\infty$ such that 
\begin{equation}\label{lim-lim}\!\!\lim_{n\rightarrow+\infty}\!\!\langle -p,x(t_n)-z\rangle\!=\!
\liminf_{t\rightarrow+\infty}\langle -p,x(t)-z\rangle\in\R.\end{equation}
Using again that $x$ is bounded, there exists $\ol x\in {\cal H}$ and a subsequence $(x(t_{n_k}))$ such that $(x(t_{n_k}))_{k\geq 0}$ 
converges weakly to $\ol x$ as $k\rightarrow+\infty$. From \eqref{lim-lim} we derive 
\begin{equation}\label{lim-olx}\liminf_{t\rightarrow+\infty}\langle -p,x(t)-z\rangle=\langle -p,\ol x-z\rangle.\end{equation}
Since $\Psi$ is weak lower semicontinuous, from Lemma \ref{l-psix=0}(ii) we get 
$$0\leq\Psi(\ol x)\leq \liminf_{k\rightarrow+\infty}\Psi(x(t_{n_k}))=0,$$
hence $\ol x\in\argmin\Psi$. Combining this with $p\in N_{\argmin\Psi}(z)$ we derive 
$\langle -p,\ol x-z\rangle\geq 0$. From \eqref{lim-olx} and \eqref{liminf-p} we conclude that \eqref{liminf-p-0} is true. 
Moreover, due to $-p-\nabla\Theta(z)\in\partial\Phi(z)$, \eqref{Phi-conv2} and \eqref{Theta-conv2} we obtain 
$$(\Phi+\Theta)(x(t))\geq (\Phi+\Theta)(z)+\langle -p,x(t)-z\rangle$$
and the conclusion follows from \eqref{liminf-p-0}. 
\end{proof}

\begin{remark}\rm One can notice that the condition $\dot\beta\leq k\beta$ has not been used in the proofs of Lemma \ref{l-psix=0} 
and Lemma \ref{l-att-cz}.
\end{remark}

We come now to the main results of the paper.

\begin{theorem}\label{main-th} Assume that $(H_\Psi)$, $(H_\Theta)$, $(H_\Phi)$, $(H^3_{\lambda})$, $(H^3_{\beta})$, $(H)$ 
and $(\widetilde{H})$ hold and let $(x,v):[0,+\infty)\rightarrow{\cal H} \times {\cal H}$ be a strong stable solution of the dynamical system 
\eqref{dyn-syst}. The following statements are true: 
\begin{enumerate}
 \item[(i)] $\int_0^{+\infty}\beta(t)\Psi(x(t))dt<+\infty$; 
 \item[(ii)] $\dot x\in L^2([0,+\infty);{\cal H})$; 
 \item[(iii)] $\langle \dot x,\dot v\rangle\in L^1([0,+\infty))$; 
 \item[(iv)] $(\Phi + \Theta)(x(t))$ converges to the optimal objective value of \eqref{opt} as $t\rightarrow+\infty$; 
 \item[(v)] $\lim_{t\rightarrow+\infty}\Psi(x(t))=\lim_{t\rightarrow+\infty}\beta(t)\Psi(x(t))=0$; 
 \item[(vi)] $x(t)$ converges weakly to an optimal solution of \eqref{opt} as $t\rightarrow+\infty$.
\end{enumerate}
\end{theorem}

\begin{proof} Take an arbitrary $z\in S$. From $(\widetilde{H})$ there exists $p\in N_{\argmin\Psi}(z)$ such that 
$-p-\nabla\Theta(z)\in\partial\Phi(z)$. Consider again the functions $g_z,h_z$ defined in Lemma \ref{l1-op}. 

Notice that statement (i) has been already proved in Lemma \ref{l1-op}.

Further, consider the function $E_2:[0,+\infty)\rightarrow \R$ defined for every $t \in [0,+\infty)$ as 
$$E_2(t)=(\Phi+\Theta)(x(t)) +\beta(t)\Psi(x(t)).$$
By using Lemma \ref{diff-Phix}, relation \eqref{dyn-syst} and $(H^3_{\beta})$ we derive for almost every $t \in [0,+\infty)$
\begin{align} \dot E_2(t)= & \ \langle v(t),\dot x(t)\rangle+\langle \nabla\Theta(x(t)),\dot x(t)\rangle+
\beta(t)\langle \nabla\Psi(x(t)),\dot x(t)\rangle+\dot\beta(t)\Psi(x(t))\nonumber\\
= & \ \langle v(t)+\nabla\Theta(x(t))+\beta(t)\nabla\Psi(x(t)),\dot x(t)\rangle+\dot\beta(t)\Psi(x(t))\nonumber\\
= & \ \langle -\lambda(t)\dot x(t)-\dot v(t),\dot x(t)\rangle+\dot\beta(t)\Psi(x(t))\nonumber\\
\leq & \ -\lambda(t)\|\dot x(t)\|^2-\langle \dot x(t),\dot v(t)\rangle + k\beta(t)\Psi(x(t))\label{diff-e2}.
\end{align}
Since $E_2$ is bounded from below, a simple integration procedure in \eqref{diff-e2} combined with (i), Lemma \ref{diff-Phix}(i)
and Lemma \ref{fejer-cont1} yields 
\begin{equation}\label{exist-lim-e2} \exists\lim_{t\rightarrow+\infty}E_2(t)\in\R,
\end{equation}
$$\int_0^{+\infty}\lambda(t)\|\dot x(t)\|^2dt<+\infty$$
and $$\int_0^{+\infty}\langle \dot x(t),\dot v(t)\rangle dt<+\infty,$$
which is statement (iii). Statement (ii) follows by taking into account that $\liminf_{t\rightarrow+\infty} \lambda(t)>0.$

Further, since $\beta(t)\Psi(x(t))\geq 0$, from \eqref{liminf} and Lemma \ref{l-att-cz} we get that 
\begin{equation}\label{liminf=0}\liminf_{t\rightarrow+\infty}\Big((\Phi +\Theta)(x(t))-(\Phi +\Theta)(z)+\beta(t)\Psi(x(t))\Big)= 0.\end{equation}
Taking into account the definition of $E_2$ and the fact that $\lim_{t\rightarrow+\infty}E_2(t)\in\R$, we conclude that 
\begin{equation}\label{lim-e2} \lim_{t\rightarrow+\infty}E_2(t)=(\Phi +\Theta)(z).
\end{equation}

Further, we have 
$$\limsup_{t\rightarrow+\infty}(\Phi +\Theta)(x(t))\leq\limsup_{t\rightarrow+\infty}\Big((\Phi +\Theta)(x(t))+\beta(t)\Psi(x(t))\Big)
=\lim_{t\rightarrow+\infty}E_2(t)=(\Phi +\Theta)(z),$$
which combined with Lemma \ref{l-att-cz} yields 
\begin{equation}\label{lim-phi} \lim_{t\rightarrow+\infty}(\Phi +\Theta)(x(t))=(\Phi +\Theta)(z),
\end{equation}
hence (iv) holds. 

The statement (v) is a consequence of Lemma \ref{l-psix=0}(ii), \eqref{lim-e2}, \eqref{lim-phi} and the definition of $E_2$. 

In order to prove statement (vi), we will make use of the Opial Lemma \ref{opial-cont}. From \eqref{ineq3} we have for almost every
$t \in [0,+\infty)$
$$\frac{\dot\lambda(t)}{2}\|x(t)-z\|^2+\lambda(t)\langle \dot x(t),x(t)-z\rangle+\frac{d}{dt}g_z(t)\leq 
\beta(t)\left(\Psi^*\left(\frac{p}{\beta(t)}\right)-\sigma_{\argmin \Psi}\left(\frac{p}{\beta(t)}\right)\right),$$
hence \begin{align}\frac{d}{dt}g_z(t)\leq & \  
\beta(t)\left(\Psi^*\left(\frac{p}{\beta(t)}\right)-\sigma_{\argmin \Psi}\left(\frac{p}{\beta(t)}\right)\right)
-\frac{\dot\lambda(t)}{2}\|x(t)-z\|^2+\lambda(t)\|\dot x(t)\|\cdot\|x(t)-z\|\nonumber\\
= & \ G_1(t)+G_2(t)\label{dgz},\end{align}
where $$G_1(t)=\beta(t)\left(\Psi^*\left(\frac{p}{\beta(t)}\right)-\sigma_{\argmin \Psi}\left(\frac{p}{\beta(t)}\right)\right)
-\frac{\dot\lambda(t)}{2}\|x(t)-z\|^2 $$
and $$G_2(t)= \lambda(t)\|\dot x(t)\|\cdot\|x(t)-z\|.$$
Now using that $x$ is bounded, from (ii) and $(H)$ we derive that $$G_1\in L^1([0,+\infty))$$ 
and $$G_2\in L^2([0,+\infty)).$$
From \eqref{dgz}, a direct application of Lemma \ref{fejer-cont2} and Lemma \ref{l1-op}(viii) yields 
$$\lim_{t\rightarrow+\infty}g_z(t)=0.$$
By combining this with Lemma \ref{l1-op}(i) and the fact that $\lim_{t\rightarrow+\infty}\lambda(t)>0$, we conclude that 
there exists $\lim_{t\rightarrow+\infty}\|x(t)-z\|\in\R$. Since $z\in S$ has been chose arbitrary, the first condition 
of the Opial Lemma is fulfilled. 

Let $(t_n)_{n\in\N}$ be a sequence of positive numbers such that $\lim_{n\rightarrow+\infty}t_n=+\infty$ and $x(t_n)$ 
converges weakly to $x_{\infty}$ as $n\rightarrow+\infty$. By using the weak lower semicontinuity of 
$\Psi$ and Lemma \ref{l-psix=0}(ii) we obtain 
$$0\leq \Psi(x_{\infty})\leq\liminf_{n\rightarrow+\infty}\Psi(x(t_n))=0,$$
hence $x_{\infty}\in \argmin\Psi$. Moreover, the weak lower semicontinuity of 
$\Phi + \Theta$ and \eqref{lim-phi} yield
$$(\Phi + \Theta)(x_{\infty})\leq\liminf_{n\rightarrow+\infty}(\Phi + \Theta)(x(t_n))=(\Phi + \Theta)(z),$$
thus $x_{\infty}\in S$.  
\end{proof}

We show in the following that if the objective function of \eqref{opt} is strongly convex, then the trajectory $x(\cdot)$
generated by \eqref{dyn-syst} converges strongly to the unique optimal solution of \eqref{opt}. 

\begin{theorem}\label{str-conv}  Assume that $(H_\Psi)$, $(H_\Theta)$, $(H_\Phi)$, $(H^3_{\lambda})$, $(H^3_{\beta})$, $(H)$ 
and $(\widetilde{H})$ hold and let $(x,v):[0,+\infty)\rightarrow{\cal H} \times {\cal H}$ be a strong stable solution of the dynamical system 
\eqref{dyn-syst}. If $\Phi+\Theta$ is strongly convex, then $x(t)$ converges strongly 
to the unique optimal solution of \eqref{opt} as $t\rightarrow+\infty$. 
\end{theorem}

\begin{proof} Let $\gamma>0$ be such that $\Phi+\Theta$ is $\gamma$-strongly 
convex. It is a well-known fact that in case the optimization problem 
\eqref{opt} has a unique optimal solution, which we denote by $z$. From $(\widetilde{H})$ there exists $p\in N_{\argmin\Psi}(z)$ 
such that $-p-\nabla\Theta(z)\in\partial\Phi(z)$. Consider again the functions $g_z,h_z$ defined in Lemma \ref{l1-op}.

By combining \eqref{ineq-conj} with the stronger inequality
\begin{equation}\label{phi-str-conv}(\Phi+\Theta)(x(t))-(\Phi+\Theta)(z)\geq 
\langle -p,x(t)-z\rangle+\frac{\gamma}{2}\|x(t)-z\|^2 \ \forall t \in [0,+\infty),\end{equation}
we obtain this time (see the proof of Lemma \ref{l1-op}) for almost every $t \in [0,+\infty)$
\begin{align} &  \frac{d}{dt}\left(\frac{\lambda(t)}{2}\|x(t)-z\|^2+g_z(t)\right)+\frac{\gamma}{2}\|x(t)-z\|^2
+\beta(t)\left(-\Psi^*\left(\frac{p}{\beta(t)}\right)+\sigma_{\argmin \Psi}\left(\frac{p}{\beta(t)}\right)\right)\nonumber\\
\leq & \ \frac{d}{dt}\left(\frac{\lambda(t)}{2}\|x(t)-z\|^2+g_z(t)\right)+\frac{\gamma}{2}\|x(t)-z\|^2+\beta(t)\Psi(x(t))+\langle -p,x(t)-z\rangle\nonumber\\
\leq & \ \frac{d}{dt}\left(\frac{\lambda(t)}{2}\|x(t)-z\|^2+g_z(t)\right)
+(\Phi+\Theta)(x(t))-(\Phi+\Theta)(z)+\beta(t)\Psi(x(t))\nonumber\\
\leq & \ 0.\label{ineq3-str-conv}
\end{align}

Taking into account $(H)$, by integration of the above inequality we obtain  
$$\int_0^{+\infty}\|x(t)-z\|^2dt<+\infty.$$ Since according to the proof of Theorem \ref{main-th}, 
$\lim_{t\rightarrow+\infty}\|x(t)-z\|$ exists, we conclude that 
$\|x(t)-z\|$ converges to $0$ as $t\rightarrow+\infty$ and the proof is complete.  
\end{proof}

\begin{remark}\label{rem13} 
The results presented in this paper remain true even if the assumed growth condition is satisfied starting
with a $t_0\geq 0$, that is, if there exists $t_0\geq 0$ such that 
$$0\leq\dot\beta(t)\leq k\beta(t) \mbox{ for almost every }  t \in [t_0,+\infty).$$
\end{remark}


\begin{thebibliography}{99}

\bibitem{abbas} B. Abbas, {\it An asymptotic viscosity selection result for the regularized Newton dynamic}, 
arXiv:1504.07793v1, 2015

\bibitem{abbas-att-sv} B. Abbas, H. Attouch, B.F. Svaiter, {\it Newton-like dynamics and forward-backward methods for 
structured monotone inclusions in Hilbert spaces}, Journal of Optimization Theory and its Applications 161(2), 331--360, 2014

\bibitem{alv-att-bolte-red} F. Alvarez, H. Attouch, J. Bolte, P. Redont, {\it A second-order gradient-like dissipative dynamical system 
with Hessian-driven damping. Application to optimization and mechanics}, Journal de Math\'{e}matiques Pures et Appliqu\'{e}es 81(8), 747--779, 2002

\bibitem{att-alv} H. Attouch, F. Alvarez, {\it The heavy ball with friction dynamical system for convex constrained 
minimization problems}, in: Optimization (Namur, 1998), 25--35, in: Lecture Notes in Economics and Mathematical Systems 481, Springer, Berlin, 2000

\bibitem{att-cab-cz} H. Attouch, A. Cabot, M.-O. Czarnecki, {\it Asymptotic behavior of nonautonomous monotone and 
subgradient evolution equations},  arXiv:1601.00767, 2016

\bibitem{att-cza-02} H. Attouch, M.-O. Czarnecki, {\it Asymptotic control and stabilization of nonlinear oscillators with 
non-isolated equilibria}, Journal of Differential Equations 179(1), 278--310, 2002

\bibitem{att-cza-10} H. Attouch, M.-O. Czarnecki, {\it Asymptotic behavior of coupled dynamical systems with multiscale aspects}, 
Journal of Differential Equations 248(6), 1315--1344, 2010

\bibitem{att-cza-16} H. Attouch, M.-O. Czarnecki, {\it Asymptotic behavior of gradient-like dynamical systems involving 
inertia and multiscale aspects},  arXiv:1602.00232, 2016

\bibitem{att-cza-peyp-p} H. Attouch, M.-O. Czarnecki, J. Peypouquet, {\it Prox-penalization and splitting methods for 
constrained variational problems}, SIAM Journal on Optimization 21(1), 149-173, 2011

\bibitem{att-cza-peyp-c} H. Attouch, M.-O. Czarnecki, J. Peypouquet, {\it Coupling forward-backward with penalty 
schemes and parallel splitting for constrained variational inequalities}, SIAM Journal on Optimization 21(4), 1251-1274, 2011

\bibitem{att-g-r} H. Attouch, X. Goudou, P. Redont, {\it The heavy ball with friction method. I. The continuous dynamical system: 
global exploration of the local minima of a real-valued function by asymptotic analysis of a dissipative dynamical system}, 
Communications in Contemporary Mathematics 2(1), 1--34, 2000

\bibitem{att-maing} H. Attouch, P.-E. Maing\'{e}, {\it Asymptotic behavior of second-order dissipative evolution equations combining 
potential with non-potential effects}, ESAIM. Control, Optimisation and Calculus of Variations 17(3), 836--857, 2011

\bibitem{att-sv2011} H. Attouch, B.F. Svaiter, {\it A continuous dynamical Newton-like approach to solving monotone inclusions}, 
SIAM Journal on Control and Optimization 49(2), 574--598, 2011

\bibitem{bauschke-book} H.H. Bauschke, P.L. Combettes, {\it Convex Analysis and Monotone Operator Theory in Hilbert Spaces}, CMS Books in Mathematics, Springer, New York, 2011

\bibitem{banert-bot-pen} S. Banert, R.I. Bo\c t, {\it Backward penalty schemes for monotone inclusion problems}, 
Journal of Optimization Theory and Applications 166(3), 930--948, 2015

\bibitem{b-hab} R.I. Bo\c t, {\it Conjugate Duality in Convex Optimization}, Lecture Notes in Economics and Mathematical Systems, 
Vol. 637, Springer, Berlin Heidelberg, 2010

\bibitem{b-c-dyn-sec-ord-pen} R.I. Bo\c t, E.R. Csetnek, {\it Second order dynamical systems associated to variational inequalities}, 
to appear in Applicable Analysis, arXiv:1512.04702v3, 2016

\bibitem{b-c-dyn-KM} R.I. Bo\c t, E.R. Csetnek, {\it A dynamical system associated with the fixed points set of a 
nonexpansive operator}, Journal of Dynamics and Differential Equations, DOI: 10.1007/s10884-015-9438-x, 2015

\bibitem{b-c-dyn-pen} R.I. Bo\c t, E.R. Csetnek, {\it Approaching the solving of constrained variational inequalities via penalty 
term-based dynamical systems}, Journal of Mathematical Analysis and Applications 435(2), 1688-1700, 2016

\bibitem{b-c-penalty-svva} R.I. Bo\c t, E.R. Csetnek, {\it Forward-backward and Tseng's type penalty schemes for monotone inclusion problems}, 
Set-Valued and Variational Analysis 22, 313--331, 2014

\bibitem{b-c-penalty-vjm} R.I. Bo\c t, E.R. Csetnek, {\it A Tseng's type penalty scheme for solving inclusion problems involving 
linearly composed and parallel-sum type monotone operators}, Vietnam Journal of Mathematics 42(4), 451--465, 2014

\bibitem{b-c-penalty-inertial} R.I. Bo\c t, E.R. Csetnek, {\it Penalty schemes with inertial effects for monotone inclusion problems}, 
arXiv:1512.04428, 2015

\bibitem{b-c-dyn-sec-ord} R.I. Bo\c t, E.R. Csetnek, {\it Second order forward-backward dynamical systems for
monotone inclusion problems}, arXiv:1503.04652, 2015

\bibitem{brezis} H. Br\'{e}zis, {\it Op\'{e}rateurs maximaux monotones et semi-groupes de contractions dans les espaces de Hilbert}, 
North-Holland Mathematics Studies No. 5, Notas de Matem\'{a}tica (50), North-Holland/Elsevier, New York, 1973

\bibitem{haraux} A. Haraux, {\it Syst\`{e}mes Dynamiques Dissipatifs et Applications},   
Recherches en Math\'{e}- matiques Appliqu\'{e}ées 17, Masson, Paris, 1991

\bibitem{noun-peyp} N. Noun, J. Peypouquet, {\it Forward-backward penalty scheme for constrained convex minimization 
without inf-compactness}, Journal of Optimization Theory and Applications, 158(3), 787--795, 2013

\bibitem{peyp-12} J. Peypouquet, {\it Coupling the gradient method with a general exterior penalization scheme for 
convex minimization}, Journal of Optimizaton  Theory and Applications 153(1), 123-138, 2012

\bibitem{sontag} E.D. Sontag, {\it Mathematical control theory. Deterministic finite-dimensional systems}, Second edition, 
Texts in Applied Mathematics 6, Springer-Verlag, New York, 1998

\bibitem{Zal-carte} C. Z\u alinescu, {\it Convex Analysis in General Vector Spaces}, World Scientific, Singapore, 2002

\end{thebibliography}
\end{document}